\documentclass[12pt]{article}

\usepackage[english]{babel}
\usepackage{amssymb,amsmath}
\usepackage{hyperref}

\usepackage{graphicx}
\usepackage{emptypage}

\newtheorem{theorem}{Theorem }[section]

\newtheorem{proposition}[theorem]{Proposition}
\newtheorem{definition}[theorem]{Definition}

\newtheorem{conjecture}[theorem]{Conjecture}
\newcommand{\proof}{\noindent\textbf{Proof. }}
\newcommand{\qed}{\hspace*{\fill}$\Box$}

\title{ Results on partial geometries with an abelian Singer group of rigid type}

\author{Stefaan De Winter\footnote{Department of Mathematical Sciences, Michigan Technological University}
  \and Ellen Kamischke\footnote{Department of Mathematical Sciences, Michigan Technological University} \and Erik Neubert\footnote{Department of Mathematical Sciences, Michigan Technological University}        \and   Zeying Wang\footnote{Department of Mathematical Sciences, Michigan Technological University} }
\date{}

\begin{document}

\maketitle

\begin{abstract}

A  partial geometry $S$ admitting an abelian Singer group $G$ is called of rigid type if all lines of $S$ have a trivial stabilizer in $G$. In this paper, we show that if a partial geometry of rigid type has fewer than $1000000$ points it must be the Van Lint-Schrijver geometry or be a hypothetical geometry with 1024 or 4096 or 194481 points, which provides evidence that partial geometries of rigid type are very rare. Along the way we also exclude an infinite set of parameters that originally seemed very promising for the construction of partial geometries of rigid type (as it contains the Van Lint-Schrijver parameters as its smallest case and one of the other cases we cannot exclude as the second member of this parameter family). We end the paper with a conjecture on this type of geometries.

\end{abstract}

\hspace*{1cm} {\bf Key words:} Partial geometry, Partial difference set, Singer group

\section{Introduction}

Ever since James Singer in his seminal paper \cite{Singer} showed that every finite Desarguesian projective space admits a cyclic sharply transitive group of automorphisms, finite geometries admitting a (not necessarily cyclic) sharply transitive group of automorphisms have attracted significant interest. Such automorphism groups are now commonly called {\it Singer groups} in honor of J. Singer.  A lot of research in this realm  has focused on projective planes, where the major conjecture is that ``{\it if a finite projective plane admits a cyclic Singer group it must be Desarguesian}''. Despite much research and some progress the conjecture is still wide open, see for example \cite{Schmidt}. Generalized quadrangles are another type of geometries that have been studied in this context. We refer to \cite{Bamberg, SDWpgst1, Tao, Ghinelli, Yoshiara} among others. Partial geometries can in some sense be thought of as sitting between projective planes and generalized quadrangles. They are first mentioned in the context of abelian Singer groups in \cite{MA94b} and a more detailed study appeared in \cite{SDWpgst2} and \cite{Leung2}. Partial geometries with nonabelian Singer groups are discussed in \cite{Swartz}. In this paper we obtain further evidence that a certain type of partial geometries with an abelian Singer group is very rare, and we formulate a conjecture on this type of geometries. 

\subsection{Definitions}

We start by providing the necessary definitions and some previous results.

\begin{definition} \label{pgdef}
A proper partial geometry, denoted $pg(s, t, \alpha)$, is a finite point-line geometry such that each line contains $s + 1$ points, each point is contained in $t+1$ lines, each pair of lines intersects in at most one point, and such that 
\begin{itemize}
\item given any point $x$ and any line $L$ not containing $x$, there are exactly $\alpha$ lines through $x$ that intersect $L$, with $0<\alpha<min(s,t)$.
\end{itemize}
\end{definition}

Partial geometries were introduced by Bose \cite{BosePG} in order to provide a setting and generalization for known characterization theorems for strongly regular graphs (it is readily checked that the point graph of a partial geometry is a strongly regular graph). Partial geometries with $\alpha=1$ are the so-called generalized quadrangles. Generalized quadrangles with an abelian Singer group are well-understood, they are the ones arising from (generalized) hyperovals; see \cite{SDWpgst1}. The natural language to study partial geometries admitting a sharply transitive abelian group of automorphisms is the language of partial difference sets. Partial difference sets were introduced by Ma in \cite{MA84}.

Let $G$ be a finite group of order $v$, and let $D\subseteq G$ be a subset of size $k$. We say $D$ is a $(v,k,\lambda,\mu)$-{\it partial difference set} (PDS) in $G$ if the expressions $gh^{-1}$, $g$, $h\in D$, $g\neq h$, represent each non-identity element in $D$ exactly $\lambda$ times, and each non-identity element of $G$ not in $D$ exactly $\mu$ times. If we further assume that $D^{(-1)}=D$  (where $D^{(s)}=\{g^s:g\in D\}$) and  $e \notin D$ (where $e$ is the identity element of $G$), then $D$ is called a {\it regular} partial difference set.  A regular PDS is called {\it trivial} if $D\cup\{e\}$ or $G\setminus D$ is a subgroup of $G$. The condition that $D$ be regular is not a very restrictive one, as  $D^{(-1)}=D$ is automatically fulfilled whenever $\lambda\neq\mu$, and $D$ is a PDS if and only if $D\cup\{e\}$ is a PDS.  In this paper we will also assume all groups are abelian. Throughout this paper we will use the following standard notations: $\beta=\lambda-\mu$ and  $\Delta=\beta^2+4(k-\mu)$.  For more information on partial difference sets, we refer the reader to a survey of Ma \cite{MA94b}. \\

The {\it Cayley graph over $G$ with connection set $D$}, denoted by Cay($G$, $D$), is the graph with the elements of $G$ as vertices, and in which two vertices $g$ and $h$ are adjacent if and only if $gh^{-1}$ belongs to $D$.  When the connection set $D$ is a regular partial difference set in $G$, Cay($G$, $D$) is a strongly regular graph.  Conversely, given a strongly regular graph admitting a sharply transitive abelian group of automorphisms $G$, it is easy to check that one obtains a regular partial difference in $G$ as follows: identify one of the vertices $v$ with the identity of $G$ and any other vertex $w$ with the unique element $g$ of $G$ such that $v^g=w$. The neighborhood of $v$ now provides a regular partial difference set in $G$. 

Suppose we have a $pg(s,t,\alpha)$ admitting an abelian automorphism group $G$ acting sharply transitively on the points and we identify the points with the elements of $G$ as above. Let $L_0$, $L_1$, $\cdots$, $L_t$  (seen as sets of element of $G$) be the $t+1$ lines through $e$, the identity element in $G$. Then the following result follows immediately:

\begin{theorem}  [\cite{MA94b}]
Suppose we have a $pg(s,t,\alpha)$ with an automorphism group $G$ acting sharply transitively on the points, with the notation above, $D=(\bigcup_{i=0}^t L_i)\setminus\{e\}$ is a regular PDS with parameters
$$(v,k,\lambda,\mu)=((s+1)(\alpha+st))/{\alpha}, \,s(t+1),\,s+t(\alpha-1)-1,\,\alpha(t-1)).$$
\end{theorem}

In \cite{SDWpgst2} the following observation was made for partial geometries admitting an abelian Singer group $G$:  the stabilizer of any line $L$ in $G$ is either trivial or has size $|L|=s+1$. Based on this pairs $(\mathcal{S},G)$ of partial geometries $\mathcal{S}$ admitting an abelian Singer group $G$ were divided into three disjoint classes:
\begin{itemize}
\item {\it spread type}: all lines have a stabilizer in $G$ of size $s+1$;
\item {\it mixed type}: some but not all lines have a stabilizer in $G$ of size $s+1$;
\item {\it rigid type}: all lines have a trivial stabilizer in $G$.
\end{itemize}

In \cite{SDWpgst2} it was shown that partial geometries of spread type are in some sense well-understood: they arise through linear representation of a so-called PG-regulus. The known examples are the partial geometries arising from maximal arcs and one other example constructed in \cite{FDCperp}. In \cite{SDWpgst2} it was also conjectured that partial geometries of mixed type with $\alpha=2$ do not exist and this conjecture was largely proven in \cite{Leung2}. The case $\alpha\geq3$ is wide open, but no examples are known. Finally, in \cite{SDWpgst2} it was shown that a partial geometry of rigid type with $\alpha=2$ must be the pg$(5,5,2)$ of Van Lint and Schrijver (see  \cite{VLS} for the original construction of this geometry using cyclotomy). In this paper we want to provide evidence that partial geometries of rigid type are very rare by showing that for $\alpha\leq1000$  at most five hypothetical parameter sets can give rise to such geometry, including the above mentioned Van Lint - Schrijver geometry. Along the way we also exclude an infinite set of parameters that originally seemed very promising for the construction of partial geometries of rigid type (as it contains the Van Lint - Schrijver parameters as its smallest case and one of the other cases we cannot exclude as the second member of this parameter family). As a consequence it follows that if a partial geometry of rigid type has fewer than $1,000,000$ points it must be the Van Lint-Schrijver geometry or have $1024$ or $4096$ or $194,481$ points. Partial results in this direction were previously obtained in the theses of the second and third author \cite{ElThesis, Neubert}.

Throughout the rest of this paper $\mathcal{S}$ will be a proper pg$(s,t,\alpha)$, $\alpha\geq2$, admitting the abelian group $G$ as a Singer group and such that the pair $(\mathcal{S},G)$ is of rigid type.

\section{Some necessary conditions}

In this section we derive some basic necessary conditions on the parameters of partial geometries of rigid type. First, we cite a theorem from \cite{SDWpgst2}. 

\begin{theorem}[\cite{SDWpgst2}]\label{Bensontype}
Let $\mathcal{S}$ be a partial geometry $pg(s, t, \alpha)$ and let $\theta$ be any automorphism of $\mathcal{S}$.
Let $f$ be the number of fixed points of $\mathcal{S}$ under $\theta$ and $g$ be the number of points $p$ of $\mathcal{S}$
for which $p$ is collinear with $p^\theta$ , where $p^\theta$ denotes the image of $p$ under $\theta$.
Then, $(1+t) f +g \equiv (1+s)(1+t) \pmod{s+t - \alpha +1}.$
\end{theorem}

We are now ready to prove the following:

\begin{theorem}\label{necessary} Let $\mathcal{S}$ be a pg$(s,t,\alpha)$ such that the pair $(\mathcal{S},G)$ is of rigid type. Then
\begin{itemize}
\item[{\rm a)}] $(s+1)(t+1) \equiv 0 \pmod{s+t-\alpha+1}$;
\item[{\rm b)}] $v=(s+1)(st+\alpha)/\alpha \equiv 0 \pmod{s+t-\alpha+1}$;
\item[{\rm c)}] $t+1=x(s+1)$ for some positive integer $x$;
\item[{\rm d)}] $x=\frac{c(s-\alpha)}{(\alpha+1)^2-c(s+1)}$ for some positive integer $c$;
\item[{\rm e)}] $s < (\alpha+1)^2 -1$.
\end{itemize}
\end{theorem}

\proof Clearly no non-identity element of $G$ has a fixed point in its action on $\mathcal{S}$. Let $g$ be an element of $G$ which maps some point $p_0$ to a non-collinear point (such $g$ obviously exists). Then we claim that $g$ maps any other point $p_1$ to a non-collinear point. Let $h$ be the unique element of $G$ that maps $p_0$ to $p_1$. As $h$ is an automorphism it preserves (non)-collinearity. Hence $p_1=p_0^h$ is not collinear with $p_0^{gh}=p_0^{hg}=p_1^g$ (note that we need $G$ to be abelian here). Hence $g$ indeed maps every point to a non-collinear point. So we obtain $f=g=0$ in Theorem \ref{Bensontype}. This implies that $(s+1)(t+1) \equiv 0 \pmod{s+t-\alpha+1}$.

\medskip

In a similar way there exists an element of $G$ that maps every point of $\mathcal{S}$ to a collinear point, so we obtain $f=0$ and $g=v=(s+1)(st+\alpha)/\alpha$ in Theorem \ref{Bensontype}. Combining with part a) it follows that $v\equiv 0 \pmod{s+t-\alpha+1}$.
\medskip

Pick any point $p_0$ and any line $L$ containing $p_0$ from $\mathcal{S}$. Assume that $p_0$, $p_1$, $\cdots$, $p_s$ are the $s+1$ points on $L$. Since $G$ acts sharply transitively on $S$, for any two collinear points $p_0$ and $p_i$, there exists a unique $g_i\in G$ such that $p_i^{g_i}=p_0$. Since the stabilizer in $G$ of $L$ is trivial, $\{L, L^{g_1}, L^{g_2}, \cdots,L^{g_s}\}$ is a set of $s+1$ distinct lines through the point $p_0$. If these are all the lines through $p_0$ we are done, if not, simply repeat the argument for one of the remaining lines through $p_0$.  As there are $t+1$ lines passing through the point $p_0$, it follows that $t+1=x(s+1)$ for some positive integer $x$.
\medskip

Combining $(s+1)(t+1) \equiv 0 \pmod{s+t-\alpha+1}$ and $t+1 = x(s+1)$ we obtain $x(s+1)^2\equiv 0 \pmod{s+t-\alpha+1}$. Replacing $t$ by $t=x(s+1)-1$ gives:
\begin{center}
$x(s+1)^2\equiv 0 $ (mod $(x+1)(s+1)-(\alpha+1)).$
\end{center}
Multiplying by $(x+1)^2$ and subtracting $x(\alpha + 1)^2$:
\begin{center}
$x\big((x+1)^2(s+1)^2 - (\alpha+1)^2\big) \equiv - x(\alpha+1)^2 \pmod{(x+1)(s+1)-(\alpha+1)}$, so

$x \big((x+1)(s+1)-(\alpha+1)\big)\big((x+1)(s+1) + (\alpha+1)\big)  \equiv - x (\alpha+1)^2 \pmod{(x+1)(s+1)-(\alpha+1)}$.
\end{center}

But now the left hand side is equivalent to 0, and we have $x(\alpha+1)^2 \equiv 0 \pmod{(x+1)(s+1) - (\alpha+1)}$. Since $x(\alpha+1)^2$ is positive, there is a positive integer $c$ such that  $x(\alpha+1)^2 = c\left((x+1)(s+1) - (\alpha+1)\right)$. Solving for $x$ yields d).

Using the fact that $\alpha<s$ it follows that the denominator in the expression from d) must be strictly positive. Hence e) follows.  \qed

We now cite three further necessary conditions. The first two were obtained by Ma:

\begin{proposition}[\cite{MA84}]\label{Ma2}
No non-trivial PDS exists in
\begin{itemize}
\item[{\rm a)}] an abelian group $G$ with a cyclic Sylow-$p$-subgroup and $o(G)\neq p$;
\item[{\rm b)}] an abelian group $G$ with a Sylow-$p$-subgroup isomorphic to $\mathbb{Z}_{p^s}\times\mathbb{Z}_{p^t}$ where $s\neq t$.
\end{itemize}
\end{proposition}

\begin{proposition}[\cite{MA84}]\label{Ma3}
If $D$ is a non-trivial regular PDS in the abelian group $G$, then $v$, $\Delta$ and $v^2/\Delta$ all have the same prime divisors. 
\end{proposition}

The third appeared in \cite{Leung2}:

\begin{proposition}  [\cite{Leung2}] \label{MLS}
Suppose there exists a proper partial geometry $\mathcal{S}$=pg{\rm (}$s$,$t$, $\alpha${\rm )} with an abelian Singer group $G$. Suppose that $\mathcal{S}$ is of rigid type with $s > 2\alpha - 1$. Then every prime divisor of $\Delta$ must  divide $s + 1$.
\end{proposition}

\section{Hypothetical cases with $\alpha\leq1000$}\label{table}

The Mathematica code below searches for all parameter sets $(s,t,\alpha)$, $2\leq\alpha\leq1000$, that survive the restrictions from Theorem \ref{necessary}, Proposition \ref{Ma2} (part a), Proposition \ref{Ma3} (the part stating that $v$ and $\Delta$ have the same prime divisors) and Proposition \ref{MLS}. Note that part a) from Proposition \ref{Ma2} is equivalent to the fact that if $v$ has multiple prime factors none of these can appear with exponent $1$. Despite the simplicity these conditions imply a number of other conditions that should be satisfied, such as the integrality of eigenvalues  for strongly regular graphs. We are not aware of other necessary conditions that could speed up the search. Given that for a partial geometry of rigid type we have that $\alpha\leq s\leq t$ it follows that $v\geq (\alpha+1)^2$. Hence the exhaustive search for $2\leq\alpha\leq1000$ implies an exhaustive search for hypothetical partial geometries of rigid type on fewer that $1,000,000$ points. Note however that also many larger values of $v$  are partially covered by the search.

{\linespread{1}
{\small
\begin{verbatim}
For[\[Alpha] = 2, \[Alpha] < 1001, \[Alpha]++,
 For[s = \[Alpha] + 1, s < (\[Alpha] + 1)^2 - 1, s++,
  For[c = 1, c < \[Alpha] + 2, c++,
   If[(\[Alpha] + 1)^2 - c*s - c > 0,
    x = (c*(s - \[Alpha]))/((\[Alpha] + 1)^2 - c*s - c);
    If[x == Ceiling[x], t = (s + 1)*x - 1;
     If[Mod[(s + 1)*(t + 1), s + t - \[Alpha] + 1] == 0,
      v = (s + 1)*(s*t + \[Alpha])/\[Alpha];
      \[Delta] = s + t - \[Alpha] + 1;
      If[Mod[v, \[Delta]] == 0,
       A = FactorInteger[v];
       B = FactorInteger[\[Delta]];
       If[Min[Last /@ A] > 1,
        If[SubsetQ[First /@ B, First /@ A],
         If[s <= 2*\[Alpha] || SubsetQ[First /@ FactorInteger[s + 1], First /@ B], 
          k = s*(t + 1);
          \[Lambda] = s - 1 + t*(\[Alpha] - 1);
          \[Mu] = \[Alpha]*(t + 1);
          Print[MatrixForm[{{s, t, \[Alpha], x, v, k, \[Lambda], \[Mu]}}]];
          If[Length[First /@ A] == 1, Print[v = Superscript @@@ A], 
           Print[v = CenterDot @@ (Superscript @@@ A)]]]]]]]]]]]]
Print["done"]
\end{verbatim}}}

This code produces the list of cases shown in Table \ref{detabel}.

\begin{table}[!htbp]
\centering
\caption{Parameters of hypothetical partial geometries of rigid type for $\alpha \leq 1000$.}\label{detabel}
\resizebox{\textwidth}{!}{%
\begin{tabular}{|c|c|c|c|c|c|c|c|c|}
\hline
Case No.	& $s$ 		& $t$ 		& $\alpha$ 	& $x$ 		& $v$                         		& $k$  		& $\lambda$ 		& $\mu$ 		\\ \hline
1 		& 5 		& 5   		& 2        	& 1   		& $3^4$				& 30   			& 9         		& 12    	\\ \hline
2		& 11 		& 23   		& 3        	& 2   		& $2^{10}$				& 264 			& 56         		& 72    	\\ \hline
3 		& 29		& 119 		& 5        	& 4   		& $2^8 \cdot 3^4$			& 3480		& 504         		& 600    	\\ \hline
4 		& 89 		& 719		& 9        	& 8   		& $2^{10} \cdot 5^4$		& 64080		& 5840     		& 6480    	\\ \hline
5 		& 39		& 39		& 15      	& 1   		& $2^{12}$				& 1560		& 584   		& 600    	\\ \hline
6 		& 305 		& 4895	& 17        	& 16   		& $2^{12} \cdot 3^8$		& 1493280		& 78624      		& 83232	\\ \hline
7 		& 1121	& 35903	& 33        	& 32   		& $2^{14} \cdot 17^4$		& 40248384		& 1150016 		& 1184832	\\ \hline
8 		& 4289	& 274559	& 65        	& 64   		& $2^{16} \cdot3^4 \cdot 11^4$ & 1177587840	& 17576064		& 17846400 	\\ \hline
9 		& 839		& 1679	& 69        	& 2   		& $2^4 \cdot5^5 \cdot 7^3$	& 1409520 		& 115010		& 115920    	\\ \hline
10		& 455		& 1367	& 95        	& 3   		& $2^{12} \cdot 3^6$		& 622440		& 128952		& 129960	\\ \hline
11		& 272		& 272		& 104        	& 1   		& $3^4 \cdot7^4$			& 74256		& 28287 		& 28392	\\ \hline
12		& 944		& 2834	& 104        	& 3   		& $3^4 \cdot 5^3 \cdot 7^4$	& 2676240		& 292845 		& 294840	\\ \hline
13		& 16769	& 2146559	& 129        	& 128		& $2^{18} \cdot5^4 \cdot13^4$ & 35995664640	& 274776320	& 276906240 \\ \hline
14		& 373		& 747		& 153        	& 2   		& $2^9 \cdot 11^3	$		& 279004		& 113916		& 114444	\\ \hline
15            &  66305           & 16974335     & 257                 &  256       & $2^{20}\cdot3^4\cdot43^4 $    & 1125483348480  & 4345496064  & 4362404352       \\ \hline
16             & 879     & 2639 & 319 & 3 & $2^{11}\cdot 5^5$ & 2320560 & 840080 & 842160  \\\hline
 17               & 7919 & 31679 & 395 & 4 & $2^6 \cdot 3^{10} \cdot 11^3$ & 250873920 & 12489444 & 12513600  \\ \hline
 18      &263681 & 135005183 & 513 & 512 & $2^{22}\cdot257^4$              & 35598301922304 & 69122917376 & 69257659392 \\ \hline
19           & 2295 & 4591 & 615 & 2 & $2^{14}\cdot  7^4$ & 10538640 & 2821168 & 2824080 \\ \hline 
20           & 1869 & 1869 & 714 & 1 &  $5^4 \cdot 11^4$		& 3495030	& 1334465		& 1335180 	\\ \hline
 21    & 1394 & 47429 & 774 & 34                  & $2^5 \cdot 5^3 \cdot 31^3$  & 66117420 & 36664010 & 36710820  \\ \hline
\end{tabular}}
\end{table}

The first case in the table corresponds to the only known existing partial geometry of rigid type: the Van Lint-Schrijver partial geometry. 
\medskip

Looking for possible patterns we observed that the following parameter family satisfies the necessary conditions used to produce the table:  $$pg(\alpha^2+\alpha-1, \alpha^3-\alpha-1, \alpha) \text{ for } \alpha = 2^n + 1 \text{ and } n=0,1,2,\hdots$$

This family seemed particularly promising given that for $n=0$ it yields the Van Lint-Schrijver parameters. For $n=1,2,3,4,5,6,7,8, 9$ the family yields case 2, 3, 4, 6, 7, 8, 13, 15 and 18 from the table respectively. However, in Section \ref{NER} we will show that unless $n=0$ or $n=1$ no member of this family can give rise to a partial geometry of rigid type.

\section{Proving nonexistence}

We describe below three general techniques for proving nonexistence of regular PDS, labeled T. 1, T. 2 and T. 3. The second and third were pioneered in a slightly less elegant form by three of the authors in \cite{SDWEKZW}.

Assume that $D$ is a nontrivial regular $(v,\,k,\,\lambda,\,\mu)$ partial difference set with $\Delta$ a perfect square in an abelian group $G$, where $G$ has at least two distinct prime divisors. Invoke the following result on ``sub-partial difference sets'' from S.L. Ma (quoted here in the form given in \cite{MA94b} (Theorem 7.1)):

\begin{proposition}[\cite{MA94b}]\label{Ma1}
Let D be a nontrivial regular $(v,\,k,\,\lambda,\,\mu)$ partial difference set  in an abelian group G. Suppose $\Delta=\delta^2$ is a  perfect square.  Let  $N$ be a subgroup of $G$ such that $\gcd(\left|N\right|, \left|G\right|/\left|N\right|)=1$ and $\left|G\right|/\left|N\right|$ is odd. Let 
$$\pi:=\gcd(|N|,\;\sqrt{\Delta}) \quad {\mbox and} \quad \theta:=\lfloor\frac{\beta+\pi}{2\pi}\rfloor.$$

Then  $D_1=N\,\cap D$ is a (not necessarily non-trivial)  regular $(v_1,\,k_1,\,\lambda_1,\,\mu_1)$ partial difference set in $N$ with
$$v_1=|N|, \; \beta_1=\lambda_1-\mu_1=\beta-2\theta \pi,  \; \Delta_1=\beta_1^2+4(k_1-\mu_1)=\pi^2,$$
and
 $$k_1=|N\cap D|=\frac{1}{2}\left[ |N| +\beta_1\pm \sqrt{(|N|+\beta_1)^2-(\Delta_1-\beta_1^2)(|N|-1)}  \right].$$

\end{proposition}

As we assumed that $G$ has at least two distinct  prime divisors we can always find a subgroup $N$ of $G$ containing a $(v_1,\,k_1,\,\lambda_1,\,\mu_1)$ partial difference set $D_1$ as in the proposition. The first non-existence result follows immediately:
\medskip

{\bf T. 1} If $(|N|+\beta_1)^2-(\Delta_1-\beta_1^2)(|N|-1)$ is not a perfect square $D$ cannot exist.
\medskip

The second non-existence argument is slightly more elaborate and depends on the so-called Local Multiplier Theorem for partial difference sets cited below:

\begin{theorem}[\cite{SDWEKZW}]\label{lmt}
Let $D$ be a regular $(v,k,\lambda,\mu)$ PDS in the abelian group $G$. Furthermore assume $\Delta$ is a perfect square.  Let $g\in G$ be an element of order $r$. Assume $s$ is a positive integer such that $\gcd(s,r)=1$. Then $g\in D$ if and only if $g^s\in D$.
\end{theorem} 

This yields:

{\bf T. 2} If $|N|=p^m$ for some prime $p$, and  $p-1$ does not divide $k_1=|N\cap D|$, then $D$ cannot exist. If $|G/N|=p^m$ for some prime $p$, and  $p-1$ does not divide $k-k_1$,  then $D$ cannot exist.
\medskip

Now note that for the difference of two elements of $D$ to belong to $N$ it is necessary and sufficient for these two elements to belong to the same coset of $N$ in $G$. Denote the distinct cosets of $N$ by $N, Nh_1,\hdots, Nh_{\frac{|G|}{|N|}-1}$ for certain elements $h_i\in G$, and set $B_i=|Nh_i\cap D|$, that is, the number of elements of the partial difference set $D$ in the $i$th coset of $N$. Double counting using the basic properties of a partial difference set yields 

$$ \sum_{i=1}^{\frac{|G|}{|N|}-1} B_i = k-k_1$$

and 

$$ \sum_{i=1}^{\frac{|G|}{|N|}-1} B_i(B_i - 1) + k_1(k_1 - 1) = k_1\cdot \lambda + (|N|-1-k_1)\cdot \mu .$$

From these two equations we can now compute the variance of the $B_i$ (up to scaling) in terms of the parameters of the PDS:

$$Var=\left( \frac{|G|}{|N|}-1\right) \sum_{i=1}^{\frac{|G|}{|N|}-1} B_i^2  -  \left(\sum_{i=1}^{\frac{|G|}{|N|}-1} B_i\right)^2,$$

which has to be non-negative. 
\medskip

{\bf T. 3} Whenever the above $Var$ is strictly negative we conclude $D$ cannot exists.

\section{Non-existence results}\label{NER}

In this section we show that most of the hypothetical parameter sets obtained in Section \ref{table} can be excluded using the techniques T.1, T.2 and T.3 described in the previous section.

\subsection{Cases excluded by T.1}

In all of the following cases one finds that the value $d=(|N|+\beta_1)^2-(\Delta_1-\beta_1^2)(|N|-1)$ from T.1 is either negative or not a perfect square. In the table below we provide the case number as well as all parameters needed to compute this value. We do not provide the value of $d$ explicitly as this is in most case a number with many digits. However we indicate with ``$<0$'' if it is strictly negative and with $\not\!\square$ if it is strictly positive but not a perfect square. Note that the value $|N|$ is sufficient to know what subgroup $N$ we consider. Also note that in some cases multiple choices for $N$ are possible but not all necessary give the desired result (this makes it complicated to efficiently incorporate T. 1 in or exhaustive search).

 \begin{tabular}{|c|c|c|c|c|c|c|c|c|c|c|}
\hline
Case	& $v$ 		& $k$ 		& $\lambda$ 	& $\mu$ 		& $|N|$                         		& $\pi$  		& $\theta$ 		& $\beta_1$ 	 & $d$	\\ \hline

9	& $2^4 \cdot 5^5 \cdot 7^3$		& 1409520		& 115010       	& 115920 		& $2^4 \cdot  5^5$	& 50	& $-9$      		& $-10$      & $\not\!\square$	\\ \hline
12              & $3^4 \cdot 5^3 \cdot 7^4$           & 2676240                 & 292845          & 294840         &  $5^3$ & $5^2$ & -40 & 5 & $<0$ \\ \hline
14              & $2^9 \cdot 11^3$      & 279004 & 113916 & 114444& $ 2^9$ & $2^3$ & -33& 0 & $\not\!\square$ \\ \hline
  16                & $2^{11} \cdot 5^5$ & 2320560 & 840080 & 842160 & $2^{11}$ & $2^7$ & -8& -32 & $<0$ \\ \hline
     17             & $2^6 \cdot 3^{10} \cdot 11^3$ & 250873920 & 12489444 & 12513600 & $2^6 \cdot 11^3 $ & $2^2 \cdot 11^2$ & -25& 44 & $<0$ \\ \hline
    21             & $2^5 \cdot 5^3 \cdot 31^3$  & 66117420 & 36664010 & 36710820 & $2^5 \cdot 31^3 $& $2 \cdot 31^2$ & -12 & -682 & $<0$ \\ 
 \hline

\end{tabular}

 \subsection{Cases excluded by T.2}
 
 We have one case in this category. The table follows the notation introduced so far with $k_1=|N\cap D|$.

\begin{tabular}{|c|c|c|c|c|c|c|c|c|c|c|}
\hline
Case	& $v$ 		& $k$ 		& $\lambda$ 	& $\mu$ 		& $|N|$                         		& $\pi$  		& $\theta$ 		& $\beta_1$ 	& $\Delta_1$ & $k_1$	\\ \hline
 
 20		& $5^4 \cdot 11^4$		& 3495030	& 1334465		& 1335180 & $11^4$ & $11^2$ & -3 & 11& $11^4$ & 8052 or 6600	\\ \hline
\end{tabular}

\medskip

Since $10$ doesn't divide $8052$, it follows that necessarily $k_1=6600$, that is, there are $6600$ elements from $N$ in $D$.  On the other hand, it then follows that $k-k_1=3,488,430$, which is not divisible by $4$, so this case cannot occur by T. 2.

 \subsection{Cases excluded by T.3}

We first consider the hypothetical infinite family described at the end of Section \ref{table}. We show that a  $pg(\alpha^2+\alpha-1, \alpha^3-\alpha-1, \alpha) \text{ for } \alpha = 2^n + 1 \text{ and } n\geq2$ of rigid type cannot exist by proving the following stronger statement:

\begin{theorem} \label{genRigidPDS_1}
Let $\alpha = 2^n + 1$ and $n \geq 2$. Then there does not exist an {\rm ($m^2$, $r(m+1)$,  $-m+r^2+3r$, $r^2+r$)}-PDS in an abelian group $G$, where $m=(\alpha+1)^2(\alpha-1)$, $r=\alpha^2-1$. 
\end{theorem}
\proof With the same notation as before we obtain the following table:

\begin{tabular}{|c|c|c|c|c|c|c|c|c|c|}
\hline
 $v$ 		& $k$ 		& $\lambda$ 	& $\mu$ 		& $|N|$                         		& $\pi$  		& $\theta$ 		& $\beta_1$ 	& $\Delta_1$ & $k_1$	\\ \hline
 
 $m^2$ 	& $r(m+1)$		& $-m+r^2+3r$ 	& $r^2+r$		& $2^{2n+4}$ &  $2^{n+2}$ & $-2^{n-2}-2^{2n-3}$ & 0 & $2^{2n+4}$ & $2^{2n+3} \pm 2^{n+1}$	\\ \hline

\end{tabular}

In the case where $k_1=2^{2n+3} + 2^{n+1}$ we obtain $$Var= -2^{4n-4}(2^n + 2)(7 \cdot 2^n - 2)(2^{3n} +15 \cdot 2^{2n} +9 \cdot 2^{n+2} + 28)<0,$$

and in the case where $k_1=2^{2n+3} - 2^{n+1}$ we obtain $$Var= -2^{2n-4}(2^n + 2)(7 \cdot 2^{2n}-3 \cdot 2^{n+1}-8)(2^{4n} +15 \cdot 2^{3n} +2^{2n+5} + 3 \cdot 2^{n+2} -16)<0.$$ This concludes the proof. \qed
\medskip

Note that $n\geq 2$ is necessary in the proof of the above theorem as for $n=0$ or $1$ the number of vertices is a prime power. The above theorem excludes existence of a rigid type partial geometry in cases 3, 4, 6, 7, 8, 13, 15 and 18 from the table in Section \ref{table}.
\medskip

One other case from the table in Section \ref{table} can be excluded using T.3. The table below provides the necessary parameters to carry out the computations:
\medskip

\begin{tabular}{|c|c|c|c|c|c|c|c|c|c|c|c|}
\hline
Case & $v$ 		& $k$ 		& $\lambda$ 	& $\mu$ 		& $|N|$                         		& $\pi$  		& $\theta$ 		& $\beta_1$ 	& $\Delta_1$ & $k_1$ & $Var$	\\ \hline
 
 10          & $2^{12}\cdot 3^6$ & 622440 & 128952 & 129960 & $ 2^{12} $& $2^6$ & -8 & 16 & $2^{12}$  & 2600 \mbox{ or} 1512 &  $<0$ \\\hline
         
\end{tabular}

\section{Some comments on the open cases} 

The four cases from Table \ref{detabel} that we could not deal with are:

\begin{table}[!htbp]
\centering
\caption{Remaining open cases.}
%resizebox{\textwidth}{!}{%
\begin{tabular}{|c|c|c|c|c|c|c|c|c|}
\hline
Case & $s$ 		& $t$ 		& $\alpha$ 	& $x$ 		& $v$                         		& $k$  		& $\lambda$ 		& $\mu$ 		\\ \hline
2		& 11 		& 23   		& 3        	& 2   		& $2^{10}$				& 264 			& 56         		& 72    	\\ \hline
5 		& 39		& 39		& 15      	& 1   		& $2^{12}$				& 1560		& 584   		& 600    	\\ \hline
11		& 272		& 272		& 104        	& 1   		& $3^4 \cdot7^4$			& 74256		& 28287 		& 28392	\\ \hline
19           & 2295 & 4591 & 615 & 2 & $2^{14}\cdot  7^4$ & 10538640 & 2821168 & 2824080 \\ \hline 

\end{tabular}
\end{table}

We first make a simple observation that also appeared in \cite{SDWpgst2} and \cite{Leung2}: {\it If $D$ is the PDS in the group $G$ arising from a partial geometry of rigid type then $D$ cannot contain an element of order two.} This is easily seen as any element of order two in $D$ would have to fix a line, contradicting the definition of rigid type.
\medskip

As an immediate consequence of this observation we see that if $(\mathcal{S},G)$ would be a pair of rigid type in any of the three open cases where $v$ is divisible by $2$ then the Sylow-2-group of $G$ cannot be elementary abelian as Proposition \ref{Ma1} shows that the Sylow-2-subgroup in these three cases must contain $k_1>0$ elements of $D$.
\medskip

We also want to note that J. Polhill in \cite{Polhill} constructed PDS with parameters $(2^{10}, 264, 56, 72)$ in an abelian group which is not elementary abelian. However, closer inspection of these PDS reveals that all of them contain an involution. Hence, these PDS cannot arise from a partial geometry of rigid type. We do want to point out a subtlety here. Although the PDS constructed by Polhill cannot arise from a partial geometry of rigid type, it is theoretically not impossible for the strongly regular graph arising from these PDS to be geometrical, or even hold a partial geometry of rigid type as there could possibly be non-isomorphic PDS yielding the same strongly regular graph.
\medskip

Using the techniques used to prove non-existence in the previous section we can derive some information on the hypothetical structure of the corresponding partial difference set in the two cases where $v$ is not a prime power. In particular one can derive that in Case 11,  the Sylow-7 subgroup of $G$ must contain $1050$ elements of $D$, whereas the Sylow-3 subgroup must contain either $24$ or $60$ elements of $D$. In the same way one can show that the Sylow-2 subgroup in Case 19 would have to contain $2064$ elements of $D$.

\section{Conclusions}

Combining the results of our search and the non-existence proofs in Section \ref{NER} we obtain:

\begin{theorem}
If $\mathcal{S}$ is a pg$(s,t,\alpha)$  of rigid type with $\alpha\leq1000$ then $\mathcal{S}$ is either the Van Lint - Schrijver partial geometry or $\mathcal{S}$ is a hypothetical geometry with parameters pg$(11,23,3)$ or pg$(39,39,15)$ or pg$(272,272,104)$ or pg$(2295,4591,615)$.
\end{theorem}

Based on our observations and attempts to construct partial geometries of rigid type in the open cases we propose the following conjecture:

\begin{conjecture}
If $\mathcal{S}$ is a pg$(s,t,\alpha)$  of rigid type then the number of points of $\mathcal{S}$ is either a power of $2$ or a power of $3$.
\end{conjecture}

Based on our observations we do not feel confident to make the stronger conjecture that the Van Lint - Schrijver geometry is the unique partial geometry of rigid type. 

\section*{Acknowledgement}

This material is based upon work done while the first author was supported by and serving at the National Science Foundation. Any opinion, findings, and conclusions or recommendations expressed in this material are those of the authors and do not necessarily reflect the views of the National Science Foundation.

\end{document}